\documentclass[12pt]{amsart}
\usepackage{amsmath}
\usepackage{latexsym}
\usepackage{amssymb}
\usepackage{amscd}

\usepackage{amsfonts}
\usepackage[all]{xy}

\newtheorem{theorem}{Theorem}

\newtheorem{conjecture}{Conjecture}

\newtheorem{proposition}{Proposition}

\newcounter{tmp}

\def\db#1{ \mathbf{D}^b({#1})}

\def\A{{\mathcal A}}
\def\B{{\mathcal B}}
\def\H{{\mathcal H}}

\def\ZZ{{\mathbb Z}}

\def\AA{{\mathbb A}}
\def\NN{{\mathbb N}}
\def\ZZ{{\mathbb Z}}
\def\QQ{{\mathbb Q}}

\def\PP{{\mathbb P}}

\def\Hom{\operatorname{Hom}}

\def\Spec{\operatorname{Spec}}

\def\wh{\widehat}

\def\CH{\operatorname{CH^{eff}}(\kk )}
\def\kk{{\mathsf k}}
\def\Sm{\operatorname{Sm}(\kk)}
\def\SmCor{\operatorname{SmCor}(\kk)}
\def\DM{\operatorname{DM^{eff}_{gm}} (\kk )}
\def\Dm{\operatorname{DM_{gm}} (\kk )}
\def\Mg{\operatorname{M}}
\def\wMg{\widetilde{\operatorname{M}}}
\def\ch{\operatorname{ch}}
\def\td{\operatorname{td}}
\def\supp{\operatorname{supp}}

\def\bR{{\mathbf R}}
\def\bL{{\mathbf L}}

\title[]{
Derived categories of coherent sheaves and motives.}

\author[]{Dmitri Orlov}

\address{Steklov Mathematical Institute RAN}

\email{orlov@mi.ras.ru}

\thanks{This work was done
with a partial financial support
from grant RFFI 05-01-01034, from the President of RF
young Russian scientists award MD-2731.2004.1,
from grant CRDF Award No RUM1-2661-MO-05, and from the Russian Science
Support Foundation.}

\date{}

\begin{document}

\maketitle
The bounded derived category of coherent sheaves $\db{X}$ is a
natural triangulated category which can be associated with  an
algebraic variety $X.$ It happens sometimes that two different
varieties have equivalent derived categories of coherent sheaves
 $\db{X}\simeq\db{Y}.$ There arises a natural question:
can one say anything about motives of $X$ and $Y$ in that case?
The first such example (see \cite{Mu}) -- abelian variety  $A$
and its dual $\wh{A}$ -- shows us that the motives
of such varieties are not necessary isomorphic. However, it seems
that the motives with rational coefficients are isomorphic

Recall a definition of the category of effective Chow motives
$\CH$ over a field $\kk.$ The category  $\CH$ can be obtained as the pseudo-abelian
envelope (i.e. as formal adding  of cokernels of all projectors)
of a category, whose objects are smooth projective schemes over
$\kk,$ and the group of morphisms from  $X$ to $Y$ is the sum
$\oplus_{X_i}A^m(X_i\times Y)$ (over all connected components
$X_i$) of the groups of cycles of codimension  $m=\dim Y$ on
$X_i\times Y$ modulo rational equivalence  (see \cite{Ma, Bl}). In
\cite{Vo} Voevodsky introduced a triangulated category of
geometric motives $\DM.$ He started with an additive category $\SmCor,$
objects of which are smooth schemes of finite type over $\kk,$ and
the group of morphisms from  $X$  to $Y$ is the free abelian group
generated by integral closed subschemes $Z\subset X\times Y$ such that
the projection on  $X$ is finite and surjective onto a connected component of
$X.$ There is a natural embedding $[-]:\Sm\to \SmCor$
of the category $\Sm$ of smooth schemes of finete type over $\kk.$
The category  $\SmCor$ is additive and one has $[X\coprod Y]=[X]\oplus[Y].$
Further, he considered the quotient of the homotopy
category $\H^b(\SmCor)$
of bounded complexes by minimal thick triangulated subcategory
$T,$ which contains all objects of the form  $[X\times\AA^1]\to
[X]$ and $[U\cap V]\to[U]\oplus[V]\to [X]$ for any open covering
$U\cup V=X.$ Triangulated category $\DM$ is defined as the pseudo-abelian
envelope of the quotient category
$\H^b(\SmCor)/T$(see~\cite{Vo,Bl}).

There exists a canonical functor  $\CH\to\DM,$ which is a full
embedding if $\kk$ admits resolution of singularities
 (\cite[4.2.6]{Vo}). Thus, it doesn't matter in which category (in  $\CH$ or in $\DM$)
 motives of smooth projective varieties are considered. Denote the motive of a variety $X$
by $\Mg(X),$ and  its motive in  the category of motives with rational coefficients
$\DM\otimes\QQ$(and in
$\CH\otimes\QQ$)   by $\Mg(X)_{\QQ}.$

\begin{conjecture}\label{G1} Let  $X$ and $Y$ be smooth projective varieties,
and let $\db{X}{\simeq}\db{Y}.$ Then the motives
$\Mg(X)_{\QQ}$ and $\Mg(Y)_{\QQ}$ are isomorphic in
$\CH\otimes\QQ$ (and in $\DM\otimes\QQ$)
\end{conjecture}

The category  $\DM$ has a tensor structure, and
$\Mg(X)\otimes\Mg(Y)=\Mg(X\times Y).$ One defines the Tate object
 $\ZZ(1)$ to be the image of the complex  $[\PP^1]\to[\Spec(\kk)]$
 placed in degree 2 and 3 and put
$M(p)=M\otimes\ZZ(1)^{\otimes p}$ for any motive $M\in\DM$ and
$p\in\NN.$ The triangulated category of geometric motives
 $\Dm$
is defined by formally inverting the functor  $-\otimes\ZZ(1)$ on
$\DM.$ The important and nontrivial fact here is the statement that the
canonical functor  $\DM\to\Dm$ is a full embedding
\cite[4.3.1]{Vo}. Therefore, we can work in the category $\Dm.$ Moreover (see \cite{Vo}), for any smooth
projective varieties $X, Y$ and for any integer  $i$ there is an isomorphism
$$
\Hom_{\Dm}(\Mg(X), \Mg(Y)(i)[2i])\cong A^{m+i}(X\times
Y),\quad\text{where}\quad m=\dim Y.
$$

\begin{conjecture} Let $X$ and $Y$ be smooth projective varieties
and   let $F: \db{X}\stackrel{}{\to}\db{Y}$ be a fully faithful
functor. Then  the motive  $\Mg(X)_{\QQ}(k)[2k]$ is a direct summand of
the motive $\Mg(Y)_{\QQ}$ for some integer $k\in \ZZ.$
\end{conjecture}

Suppose, one has a fully faithful functor  $F:\db{X}\to\db{Y}$
between derived categories of coherent sheaves of two smooth
projective varieties  $X$ and $Y$ of dimension $n$ and  $m$
respectively. Any such functor has a right adjoint  $F^*$ by
\cite{BV}, and by Theorem  2.2 from \cite{Or1} (see also
\cite[3.2.1]{Or2}) the functor  $F$ can be represented by an
object on the product $X\times Y,$ i.e. $F\cong\Phi_{\A},$ where
$\Phi_{\A}=\bR p_{2*}(p^*_1(-)\stackrel{\bL}{\otimes}\A)$ for some
$\A\in\db{X\times Y}.$ With any functor of the form
$\Phi_{\A}:{\db{X}}\to{\db{Y}}$ one can associate an element $a\in
A^*(X\times Y, \QQ)$ by the following rule
\begin{equation}\label{epsi}
a= p^*_1\sqrt{ \td_X}\cdot \ch(\A)\cdot p^*_2\sqrt{\td_Y},
\end{equation}
where $\td_X$ and $\td_Y$ are Todd classes of the varieties $X$ and $Y.$
The cycle $a$ has a mixed type. Let us consider its decomposition into components
 $a=a_0+\cdots+a_{n+m},$ where index is the codimension of a cycle on
 $X\times Y.$ Each component  $a_q$ induces a map of motives
$$
\alpha_q : \Mg(X)_{\QQ}\to \Mg(Y)_{\QQ}(q-m)[2(q-m)].
$$
Thus the total cycle  $a$ gives a map $ \alpha:
\Mg(X)_{\QQ}\stackrel{}{\to}\bigoplus_{i=-m}^{n}
\Mg(Y)_{\QQ}(i)[2i].$ Now consider the object $\B\in \db{X\times
Y}$ which represents the (left) adjoint functor  $F^*,$ i.e.
$F^*\cong \Psi_{\B},$ where $\Psi_{\B}=\bR
p_{1*}(p^*_2(-)\stackrel{\bL}{\otimes}\B).$
One attaches to the object
$\B$ a cycle $b=b_0+\cdots + b_{n+m}$ defined by the same formula
(\ref{epsi}). The cycle $b$ induces a map  $\beta:
\bigoplus_{i=-m}^{n} \Mg(Y)_{\QQ}(i)[2i] \stackrel{}{\to}
\Mg(X)_{\QQ}.$ Since the functor  $\Phi_{\A}$ is fully faithful, the composition
 $\Psi_{\B}\circ\Phi_{\A}$ is isomorphic to the identity functor.
 Applying the Riemann-Roch-Grothendieck theorem, we obtain that the composition
$$
\Mg(X)_{\QQ}\stackrel{\alpha}{\to}\bigoplus_{i=-m}^{n}
\Mg(Y)_{\QQ}(i)[2i] \stackrel{\beta}{\to} \Mg(X)_{\QQ}
$$
is the identity map, i.e. $\Mg(X)_{\QQ}$ is a direct summand of
 $\bigoplus_{i=-m}^{n} \Mg(Y)_{\QQ}(i)[2i].$
 
 Denote by $\wMg(X)_{\QQ}$ the infinite direct sum $\bigoplus_{i=-\infty}^{\infty} \Mg(X)_{\QQ}(i)[2i].$
 The total cycle  $a$ defined above gives a map $ \widetilde{\alpha}:
\wMg(X)_{\QQ}{\to}
\wMg(Y)_{\QQ}$ whose component from $\Mg(X)_{\QQ}(k)[2k]$ to $\Mg(Y)_{\QQ}(q-m+k)[2(q-m+k)]$ coincides with $\alpha_q$
for any $k.$ By the same way the cycle $b$ induces a map $\widetilde{\beta}$ from $\wMg(X)_{\QQ}$ to
$\wMg(Y)_{\QQ}.$ The above consideration give us that the composition $\widetilde{\beta}\cdot \widetilde{\alpha}$ is the identity map.
Thus we obtain the following proposition.

\begin{proposition}
Let $X$ and $Y$ be smooth projective varieties and let $F: \db{X}\stackrel{}{\to}\db{Y}$ be a fully
faithful functor.
Then the motive  $\wMg(X)_{\QQ}$ is a direct summand of the motive  $\wMg(Y)_{\QQ}.$
If, in addition, the functor  $F$ is an equivalence,
then the motives  $\wMg(X)_{\QQ}$ and $\wMg(Y)_{\QQ}$ are isomorphic.
\end{proposition}

Assume now that $\dim X=\dim Y=n$ and, moreover, suppose that the support of the
object
$A$ also has the dimension  $n.$
Therefore, $a_q=0$ when $q=0,\dots,
n-1,$ i.e. $a=a_n+\cdots +a_{2n}.$
It is easily to see that in this case $b=b_{n}+\cdots+b_{2n}$ as well.
This implies that the composition
$\beta\cdot\alpha: \Mg(X)_{\QQ}\to \Mg(X)_{\QQ},$ which is the identity,
coincides with $\beta_n\cdot\alpha_n$. Hence,
$\Mg(X)_{\QQ}$ is a direct summand of  $\Mg(Y)_{\QQ}.$
Furthermore, since the cycles  $a_n$ and $b_n$ are integral in this case we get
the same result for integral motives, i.e. the integral motive  $\Mg(X)$ is a direct summand of
the motive $\Mg(Y)$ as well. Thus, we obtain

\begin{theorem}
Let $X$ and $Y$ be smooth projective varieties of dimension
$n,$ and let $F: \db{X}\stackrel{}{\to}\db{Y}$ be a fully
faithful functor such that the dimension of the support of an object
$\A$ on $X\times Y,$  which represents $F,$ is equal to $n.$
Then the motive  $\Mg(X)$ is a direct summand of the motive  $\Mg(Y).$
If, in addition, the functor  $F$ is an equivalence,
then the motives  $\Mg(X)$ and $\Mg(Y)$ are isomorphic.
\end{theorem}

Examples of such functors are known, they come from birational
geometry (see e.g. \cite{Or2}). In these examples
one of the connected components of  $\supp(\A)$ gives a birational map
$X\dasharrow Y.$ Blow ups and antiflips induce fully faithful functors,
and flops induce equivalences.
Note that an isomorphism of motives implies an isomomorphism of
any realization (singular cohomologies, l-adic cohomologies, Hodge structures and  so on).

For arbitrary equivalence  $\Phi_{\A}: \db{X}\to\db{Y}$
the map of motives $\alpha_n:\Mg(X)_{\QQ}\to \Mg(Y)_{\QQ},$
induced by the cycle  $a_n\in A^n(X\times Y, \QQ),$ is not necessary an isomorphism (e.g. Poincare
line bundle $\mathcal{P}$ on the product of abelian variety $A$ and its dual $\wh{A}$).
However, the following conjecture, which specifies Conjecture \ref{G1}, may be true.

\begin{conjecture} Let $\A$ be an object of $\db{X\times Y},$
for which
$\Phi_{\A}:\db{X}\to\db{Y}$ is an equivalence.
Then there exist line bundles  $L$ and $M$ on $X$ and on $Y$ respectively
such that the component
$a'_{n}$ of the object $\A':= p_1^* L\otimes \A\otimes p_2^*
M$ gives an isomorphism between motives  $\Mg(X)_{\QQ}$ and
$\Mg(Y)_{\QQ}.$
\end{conjecture}

I am grateful to Yu. I. Manin for very
useful discussions.

\end{document}